\documentstyle[amssymb]{article}

\newtheorem{th}{Theorem}[section]

\newtheorem{prop}[th]{Proposition}
\newtheorem{cor}[th]{Corollary}
\newtheorem{defn}[th]{Definition}
\newenvironment{defn-new}{\begin{defn} \em}{\end{defn}}
\newtheorem{rem}[th]{Remark}
\newenvironment{rem-new}{\begin{rem} \em}{\end{rem}}
\newtheorem{ex}[th]{Example}
\newenvironment{ex-new}{\begin{ex} \em}{\end{ex}}
\newtheorem{exer}[th]{Exercise}
\newenvironment{exer-new}{\begin{exer} \em}{\end{exer}}
\newtheorem{agr}[th]{Agreement}
\newenvironment{agr-new}{\begin{agr} \em}{\end{agr}}
\newtheorem{pbm}[th]{Problem}
\newenvironment{pbm-new}{\begin{pbm} \em}{\end{pbm}}

\makeatletter \@addtoreset{equation}{section} \makeatother

\begin{document}

\begin{center}
{\Large {\bf Einstein like $(\varepsilon )$-para Sasakian manifolds}}
\bigskip

Sadik Kele\c{s}, Erol Kili\c{c}, Mukut Mani Tripathi and Selcen Y\"{u}ksel
Perkta\c{s}\bigskip \bigskip
\end{center}

\noindent {\bf Abstract.} Einstein like $(\varepsilon )$-para Sasakian
manifolds are introduced. For an $\left( \varepsilon \right) $-para Sasakian
manifold to be Einstein like, a necessary and sufficient condition in terms
of its curvature tensor is obtained. The scalar curvature of an Einstein
like $\left( \varepsilon \right) $-para Sasakian manifold is obtained and it
is shown that the scalar curvature in this case must satisfy certain
differential equation. A necessary and sufficient condition for an $\left(
\varepsilon \right) $-almost paracontact metric hypersurface of an
indefinite locally Riemannian product manifold to be $\left( \varepsilon
\right) $-para Sasakian is obtained and it is proved that the $\left(
\varepsilon \right) $-para Sasakian hypersurface of an indefinite locally
Riemannian product manifold of almost constant curvature is always Einstein
like. \medskip

 \medskip


\section{Introduction\label{sect-intro}}

In 1976, S\={a}to \cite{Sato-76} introduced an almost paracontact structure
on a differentiable manifold, which is an analogue of the almost contact
structure \cite{Sasaki-60-Tohoku,Blair-02-book} and is closely related to
almost product structure (in contrast to almost contact structure, which is
related to almost complex structure). An almost contact manifold is always
odd-dimensional but an almost paracontact manifold could be even-dimensional
as well. In 1969, Takahashi \cite{Takahashi-69-Tohoku-1} studied almost
contact manifolds equipped with associated pseudo-Riemannian metrics. The
indefinite almost contact metric manifolds and indefinite Sasakian manifolds
are also known as $\left( \varepsilon \right) $-almost contact metric
manifolds and $\left( \varepsilon \right) $-Sasakian manifolds, respectively
\cite{Bej-Dug-93,Duggal-90-IJMMS}. Also, in 1989, Matsumoto \cite{Mat-89}
replaced the structure vector field $\xi $ by $-\,\xi $ in an almost
paracontact manifold and associated a Lorentzian metric with the resulting
structure and called it a Lorentzian almost paracontact manifold. In a
Lorentzian almost paracontact manifold given by Matsumoto, the
semi-Riemannian metric has only index $1$ and the structure vector field $%
\xi $ is always timelike. Because of these circumstances, the authors in
\cite{TKYK-09-IJMMS} introduced $\left( \varepsilon \right) $-almost
paracontact structures by associating a semi-Riemannian metric, not
necessarily Lorentzian, with an almost paracontact structure, where the
structure vector field $\xi $ is spacelike or timelike according as $%
\varepsilon =1$ or $\varepsilon =-1$. \medskip

In \cite{Sharma-82}, Sharma introduced and studied Einstein like para
Sasakian manifolds. Motivated by his study, in this paper we introduce and
study Einstein like $(\varepsilon )$-almost paracontact metric manifolds.
The paper is organized as follows. Section~\ref{sect-prel} contains some
preliminaries about $\left( \varepsilon \right) $-para Sasakian manifolds.
In section~\ref{sect-eta-phi-Einst}, we give the definition of an Einstein
like $\left( \varepsilon \right) $-almost paracontact metric manifold and
give some basic properties. For an $\left( \varepsilon \right) $-para
Sasakian manifold to be Einstein like, we also find a necessary and
sufficient condition in terms of its curvature tensor. We also find the
scalar curvature of an Einstein like $\left( \varepsilon \right) $-para
Sasakian manifold and show that the scalar curvature in this case must
satisfy certain differential equation. In section~\ref{sect-eps-PS-hyp}, we
find a necessary and sufficient condition for an $\left( \varepsilon \right)
$-almost paracontact metric hypersurface of an indefinite locally Riemannian
product manifold to be $\left( \varepsilon \right) $-para Sasakian. Finally
we prove that an $\left( \varepsilon \right) $-para Sasakian hypersurface of
an indefinite locally Riemannian product manifold of almost constant
curvature is always Einstein like.

\section{Preliminaries\label{sect-prel}}

Let $M$ be an $n$-dimensional almost paracontact manifold \cite{Sato-76}
equipped with an almost paracontact structure $(\varphi ,\xi ,\eta )$
consisting of a tensor field $\varphi $ of type $(1,1)$, a vector field $\xi
$ and a $1$-form $\eta $ satisfying
\[
\varphi ^{2}=I-\eta \otimes \xi ,\quad \eta (\xi )=1,\quad \varphi \xi
=0,\quad \eta \circ \varphi =0.
\]%
By a semi-Riemannian metric \cite{ONeill-83} on a manifold $M$, we
understand a non-degenerate symmetric tensor field $g$ of type $\left(
0,2\right) $. In particular, if its index is $1$, it becomes a Lorentzian
metric \cite{Beem-Ehrlich-81}. Throughout the paper we assume that $%
X,Y,Z,U,V,W\in {\frak X}\left( M\right) $, where ${\frak X}\left( M\right) $
is the Lie algebra of vector fields in $M$, unless specifically stated
otherwise. Let $g$ be a semi-Riemannian metric with ${\rm index}(g)=\nu $ in
an $n$-dimensional almost paracontact manifold $M$ such that
\begin{equation}
g\left( \varphi X,\varphi Y\right) =g\left( X,Y\right) -\varepsilon \eta
(X)\eta \left( Y\right) ,  \label{eq-metric-1}
\end{equation}%
where $\varepsilon =\pm 1$. Then $M$ is called an $\left( \varepsilon
\right) $-almost paracontact metric manifold equipped with an $\left(
\varepsilon \right) ${\em -}almost paracontact metric structure $(\varphi
,\xi ,\eta ,g,\varepsilon )$ \cite{TKYK-09-IJMMS}. In particular, if ${\rm %
index}(g)=1$, then an $(\varepsilon )$-almost paracontact metric manifold is
a Lorentzian almost paracontact manifold. In particular, if the metric $g$
is positive definite, then an $(\varepsilon )$-almost paracontact metric
manifold is the usual almost paracontact metric manifold \cite{Sato-76}. The
equation (\ref{eq-metric-1}) is equivalent to
\begin{equation}
g\left( X,\varphi Y\right) =g\left( \varphi X,Y\right) \quad {\rm along\ with%
}\quad g\left( X,\xi \right) =\varepsilon \eta (X).  \label{eq-metric-2}
\end{equation}%
Note that $g\left( \xi ,\xi \right) =\varepsilon $, that is, the structure
vector field $\xi $ is never lightlike. An $\left( \varepsilon \right) $%
-almost paracontact metric structure $(\varphi ,\xi ,\eta ,g,\varepsilon )$
is called an $\left( \varepsilon \right) ${\em -}para Sasakian structure if
\begin{equation}
\left( \nabla _{X}\varphi \right) Y=-\,g(\varphi X,\varphi Y)\xi
-\varepsilon \eta \left( Y\right) \varphi ^{2}X,  \label{eq-eps-PS-def-1}
\end{equation}%
where $\nabla $ is the Levi-Civita connection with respect to $g$. A
manifold endowed with an $\left( \varepsilon \right) $-para Sasakian
structure is called an $\left( \varepsilon \right) $-para Sasakian manifold.
In an $\left( \varepsilon \right) $-para Sasakian manifold we have
\begin{equation}
\nabla \xi =\varepsilon \varphi ,  \label{eq-s-pcm-def}
\end{equation}%
\begin{equation}
\Phi \left( X,Y\right) \equiv g\left( \varphi X,Y\right) =\varepsilon
g\left( \nabla _{X}\xi ,Y\right) =\left( \nabla _{X}\eta \right) Y.
\label{eq-s-pcm-def1}
\end{equation}%
For more details we refer to \cite{TKYK-09-IJMMS}.

\section{Einstein like $\left( \protect\varepsilon \right) $-para Sasakian
manifolds\label{sect-eta-phi-Einst}}

We begin with the following definition analogous to Einstein like para
Sasakian manifolds \cite{Sharma-82}.

\begin{defn-new}
An $\left( \varepsilon \right) $-almost paracontact metric manifold is said
to be Einstein like if its Ricci tensor $S$ satisfies
\begin{equation}
S\left( X,Y\right) =a\,g\left( X,Y\right) +b\,g\left( \varphi X,Y\right)
+c\,\eta \left( X\right) \eta \left( Y\right)   \label{eq-(eta-phi)-Einst-1}
\end{equation}%
for some real constants $a$, $b$ and $c$.
\end{defn-new}

\begin{prop}
In an Einstein like $\left( \varepsilon \right) $-almost paracontact metric
manifold, we have
\begin{equation}
S\left( \varphi X,Y\right) =ag\left( \varphi X,Y\right) +bg\left( \varphi
X,\varphi Y\right) ,  \label{eq-(eta-phi)-Einst-3}
\end{equation}%
\begin{equation}
S\left( X,\xi \right) =\varepsilon a\eta \left( X\right) +c\eta \left(
X\right) ,  \label{eq-(eta-phi)-Einst-4}
\end{equation}%
Moreover, if the manifold is $\left( \varepsilon \right) $-para Sasakian,
then
\begin{equation}
\varepsilon a+c=1-n,  \label{eq-a-c-n}
\end{equation}%
\begin{equation}
r=na+b\,{\rm trace}(\varphi )+\varepsilon c,  \label{eq-(eta-phi)-Einst-r-1}
\end{equation}%
where $r$\ is the scalar curvature.
\end{prop}

\noindent {\bf Proof.} The equations (\ref{eq-(eta-phi)-Einst-3}) and (\ref%
{eq-(eta-phi)-Einst-4}) are obvious. In an $\left( \varepsilon \right) $%
-para Sasakian manifold, it follows that $S(X,\xi )=-(n-1)\eta (X)$, which
in view of (\ref{eq-(eta-phi)-Einst-4}) implies (\ref{eq-a-c-n}). Now, let $%
\left\{ e_{1},\ldots ,e_{n}\right\} $ be a local orthonormal frame. Then
from (\ref{eq-(eta-phi)-Einst-1}), we have
\[
r=\sum_{i=1}^{n}\left\{ \varepsilon _{i}ag\left( e_{i},e_{i}\right)
+\varepsilon _{i}bg\left( \varphi e_{i},e_{i}\right) +\varepsilon
_{i}cg\left( \xi ,e_{i}\right) g\left( \xi ,e_{i}\right) \right\} ,
\]%
which gives (\ref{eq-(eta-phi)-Einst-r-1}). $\blacksquare $

\begin{th}
In an Einstein like $\left( \varepsilon \right) $-para Sasakian manifold,
the scalar curvature $r$\ satisfies the following two differential equations
\begin{equation}
b\,\xi r-2cr=2\varepsilon \left( 1-n\right) \left( b^{2}-c^{2}-cn\right) .
\label{eq-(eta-phi)-Einst-r-4}
\end{equation}
\end{th}

\noindent {\bf Proof.} From (\ref{eq-(eta-phi)-Einst-1}), it follows that
the Ricci operator $Q$ satisfies
\begin{equation}
QX=aX+b\varphi X+\varepsilon c\,\eta \left( X\right) \xi .
\label{eq-(eta-phi)-Einst-2}
\end{equation}%
Differentiating (\ref{eq-(eta-phi)-Einst-2}), we find
\[
\left( \nabla _{Y}Q\right) X=b\left( \nabla _{Y}\varphi \right)
X+\varepsilon c\left( \nabla _{Y}\eta \right) \left( X\right) \xi
+\varepsilon c\eta \left( X\right) \nabla _{Y}\xi .
\]%
Using (\ref{eq-eps-PS-def-1}), (\ref{eq-s-pcm-def1}) and (\ref{eq-s-pcm-def}%
) in the above equation we get
\begin{eqnarray}
\left( \nabla _{Y}Q\right) X &=&-\,\varepsilon b\eta \left( X\right) Y+c\eta
\left( X\right) \varphi Y  \nonumber \\
&&-\left( bg\left( X,Y\right) -2\varepsilon b\eta \left( X\right) \eta
\left( Y\right) -\varepsilon cg\left( \varphi X,Y\right) \right) \xi .
\label{eq-del-Q-1}
\end{eqnarray}%
Now, using (\ref{eq-del-Q-1}) we have
\begin{equation}
\left( {\rm div}Q\right) X=\left\{ \varepsilon \left( 1-n\right) b+c{\rm %
\,trace}(\varphi )\right\} \eta \left( X\right) .  \label{eq-div-Q-1}
\end{equation}%
From (\ref{eq-(eta-phi)-Einst-r-1}) and (\ref{eq-a-c-n}) we get
\begin{equation}
r=b\,{\rm trace}(\varphi )-\varepsilon \left( n-1\right) \left( c+n\right)
\label{eq-(eta-phi)-Einst-r-2}
\end{equation}%
Using $Xr=2\left( {\rm div}Q\right) X$ and (\ref{eq-(eta-phi)-Einst-r-2}) in
(\ref{eq-div-Q-1}) we get (\ref{eq-(eta-phi)-Einst-r-4}). $\blacksquare $

\begin{th}
In an Einstein like $\left( \varepsilon \right) $-para Sasakian manifold, if
${\rm trace}(\varphi )$ is constant then
\begin{equation}
{\rm trace}(\varphi )=\frac{\varepsilon \left( n-1\right) b}{c}.
\label{eq-tr-phi}
\end{equation}
\end{th}

\noindent {\bf Proof. }Using $Xr=2\left( {\rm div}Q\right) X$ in (\ref%
{eq-div-Q-1}), we get
\begin{equation}
dr=2\left( \varepsilon \left( 1-n\right) b+c{\rm \,trace}(\varphi )\right)
\eta .  \label{eq-dr}
\end{equation}%
Since ${\rm trace}(\varphi )$ is constant, from (\ref{eq-(eta-phi)-Einst-r-1}%
), it follows that $r$ is constant. Hence (\ref{eq-dr}) gives (\ref%
{eq-tr-phi}). $\blacksquare $ \medskip

From now on in this section the ${\rm trace}(\varphi )$ will be assumed to
be constant.

\begin{th}
An $\left( \varepsilon \right) $-para Sasakian manifold with constant ${\rm %
trace}(\varphi )$ is Einstein like if and only if the $\left( 0,2\right) $%
-tensor field $C_{1}^{1}(\varphi R)$ is a linear combination of $g$, $\Phi $
and $\eta \otimes \eta $ formed with constant coefficients.
\end{th}

\noindent {\bf Proof. }In an $\left( \varepsilon \right) $-para Sasakian
manifold the curvature tensor $R$ satisfies \cite{TKYK-09-IJMMS}
\begin{eqnarray*}
R\left( X,Y\right) \varphi Z &=&\varphi R\left( X,Y\right) Z+\varepsilon
\Phi \left( Y,Z\right) X-\varepsilon \,\Phi \left( X,Z\right) Y-2\varepsilon
\,\Phi \left( Y,Z\right) \eta \left( X\right) \xi +2\varepsilon \,\Phi
\left( X,Z\right) \eta \left( Y\right) \xi  \\
&&-\varepsilon g\left( Y,Z\right) \varphi X+\varepsilon g\left( X,Z\right)
\varphi Y+2\eta \left( Y\right) \eta \left( Z\right) \varphi X-2\eta \left(
X\right) \eta \left( Z\right) \varphi Y.
\end{eqnarray*}%
Then we have
\begin{equation}
S\left( Y,\varphi Z\right) =C_{1}^{1}(\varphi R)\left( Y,Z\right)
+\varepsilon \left( n-2\right) \Phi \left( Y,Z\right) +\left( 2\eta \left(
Y\right) \eta \left( Z\right) -\varepsilon g\left( Y,Z\right) \right) {\rm %
trace}(\varphi ).  \label{eq-S(Y,phiZ)-1}
\end{equation}%
Since in an $\left( \varepsilon \right) $-para Sasakian manifold, it follows
that \cite{TKYK-09-IJMMS} $S\left( X,\varphi Y\right) =S\left( \varphi
X,Y\right) $, and also it can be verified that $C_{1}^{1}(\varphi R)\left(
Y,Z\right) =C_{1}^{1}(\varphi R)\left( Z,Y\right) $; therefore the equation (%
\ref{eq-S(Y,phiZ)-1}) is consistent. Now, if the manifold is Einstein like
then from (\ref{eq-(eta-phi)-Einst-3}), (\ref{eq-S(Y,phiZ)-1}) and (\ref%
{eq-tr-phi}), it follows that
\begin{equation}
C_{1}^{1}(\varphi R)=\frac{b}{c}\left( c+n-1\right) g+\left( a-\varepsilon
\left( n-2\right) \right) \Phi -\frac{\varepsilon }{c}\left( c+2b\left(
n-1\right) \right) \eta \otimes \eta ,  \label{eq-C11-phiR-1}
\end{equation}%
which shows that $C_{1}^{1}(\varphi R)$ is a linear combination of $g$, $%
\Phi $ and $\eta \otimes \eta $ formed with constant coefficients. The
converse is easy to follow. $\blacksquare $

\begin{cor}
In an Einstein like $\left( \varepsilon \right) $-para Sasakian manifold
with constant ${\rm trace}(\varphi )$, the $\left( 0,2\right) $-tensor field
$C_{1}^{1}(\varphi R)$ is parallel along the vector field $\xi $.
\end{cor}

\noindent {\bf Proof.} Since in an Einstein like $\left( \varepsilon \right)
$-para Sasakian manifold $\nabla _{\xi }\Phi =0$ and $\nabla _{\xi }\eta =0$%
, therefore from (\ref{eq-C11-phiR-1}) we conclude that $C_{1}^{1}(\varphi
R) $ is parallel along the vector field $\xi $. $\blacksquare $

\begin{th}
In an Einstein like $\left( \varepsilon \right) $-para Sasakian manifold, we
have
\begin{equation}
{\frak L}_{\xi }S=2a\varepsilon \Phi +2b\varepsilon \left( g-\varepsilon
\eta \otimes \eta \right) .  \label{eq-Lie-S}
\end{equation}
\end{th}

\noindent {\bf Proof.} In an $\left( \varepsilon \right) $-para Sasakian
manifold, we obtain
\begin{equation}
{\frak L}_{\xi }\eta =\nabla _{\xi }\eta =0,\quad {\frak L}_{\xi }\Phi
=2\varepsilon \left( g-\eta \otimes \eta \right) ,\quad {\frak L}_{\xi
}g=2\varepsilon \Phi .  \label{eq-Lie-eta}
\end{equation}%
Now, taking Lie derivative of $S$ in the direction of $\xi $ in (\ref%
{eq-(eta-phi)-Einst-1}) and using (\ref{eq-Lie-eta}), we obtain (\ref%
{eq-Lie-S}). $\blacksquare $

\begin{th}
In an Einstein like $\left( \varepsilon \right) $-para Sasakian manifold
with constant ${\rm trace}(\varphi )$, we have
\begin{equation}
{\frak L}_{\xi }\left( C_{1}^{1}(\varphi R)\right) =\frac{2\varepsilon b}{c}%
\left( c+n-1\right) \Phi +2\varepsilon \left( a-\varepsilon \left(
n-2\right) \right) \left( g-\eta \otimes \eta \right) .
\label{eq-Lie-C11-phiR}
\end{equation}
\end{th}

\noindent {\bf Proof.} Taking Lie derivative of $C_{1}^{1}(\varphi R)$ in
the direction of $\xi $ in (\ref{eq-C11-phiR-1}) and using (\ref{eq-Lie-eta}%
), we get (\ref{eq-Lie-C11-phiR}). $\blacksquare $

\section{$\left( \protect\varepsilon \right) $-para Sasakian hypersurfaces
\label{sect-eps-PS-hyp}}

Let $\tilde{M}$ be a real $\left( n+1\right) $-dimensional manifold. Suppose
$\tilde{M}$ is endowed with an almost product structure $J$ and a
semi-Riemannian metric $\tilde{g}$ satisfying
\begin{equation}
\tilde{g}(JX,JY)=\tilde{g}\left( X,Y\right)  \label{eq-g(JX,JY)}
\end{equation}%
for all vector fields $X,Y$ in $\tilde{M}$. Then we say that $\tilde{M}$ is
an indefinite almost product Riemannian manifold. Moreover, if on $\tilde{M}$
we have
\begin{equation}
(\tilde{\nabla}_{X}J)Y=0  \label{eq-del-J=0}
\end{equation}%
for all $X,Y\in {\frak X}(\tilde{M})$, where $\tilde{\nabla}$ is the
Levi-Civita connection with respect to $\tilde{g}$, we say that $\tilde{M}$
is an indefinite locally Riemannian product manifold. \medskip

Now, let $M$ be an orientable non-degenerate hypersurface of $\tilde{M}$.
Suppose now that $N$ is the normal unit vector field of $M$ such that $%
\tilde{g}\left( N,N\right) =\varepsilon $ and
\begin{equation}
JN=\xi \in {\frak X}(M).  \label{eq-JN=xi}
\end{equation}%
Let
\begin{equation}
JX=\varphi X+\eta \left( X\right) N.  \label{eq-JX}
\end{equation}

\begin{prop}
The set $\left( \varphi ,\xi ,\eta ,g\right) $ is an $\left( \varepsilon
\right) $-almost paracontact metric structure, where $g$ is the induced
metric on $M$.
\end{prop}

\noindent {\bf Proof.} We have
\[
X=J^{2}X=\varphi ^{2}X+\eta \left( \varphi X\right) N+\eta \left( X\right)
\xi ,
\]%
where (\ref{eq-JX}) and (\ref{eq-JN=xi}) are used. Equating tangential and
normal parts we get $\varphi ^{2}=I-\eta \otimes \xi $ and $\eta \circ
\varphi =0$, respectively. We also have
\[
N=J^{2}N=J\xi =\varphi \xi +\eta \left( \xi \right) N,
\]%
where (\ref{eq-JX}) and (\ref{eq-JN=xi}) are used. Equating tangential and
normal parts we get $\varphi \xi =0$ and $\eta \left( \xi \right) =1$,
respectively. Finally, we have $g\left( X,Y\right) =\tilde{g}\left(
JX,JY\right) $, which in view of (\ref{eq-JX}) gives (\ref{eq-metric-1}). $%
\blacksquare $ \medskip

The Gauss and Weingarten formulas are given respectively by
\begin{equation}
\tilde{\nabla}_{X}Y=\nabla _{X}Y+\varepsilon g(AX,Y)N,  \label{eq-GF}
\end{equation}%
\begin{equation}
\tilde{\nabla}_{X}N=-AX,  \label{eq-WF}
\end{equation}%
where $\nabla $ is the Levi-Civita connection with respect to the
semi-Riemannian metric $g$ induced by $\tilde{g}$ on $M$ and $A$ is the
shape operator of $M$.

\begin{prop}
The $\left( \varepsilon \right) $-almost paracontact metric structure on $M$
satisfies
\begin{equation}
\left( \nabla _{X}\varphi \right) Y=\eta (Y)AX+\varepsilon g(AX,Y)\xi ,
\label{eq-del-phi}
\end{equation}%
\begin{equation}
\left( \nabla _{X}\eta \right) Y=-\varepsilon g(AX,\varphi Y),
\label{eq-del-eta}
\end{equation}%
\begin{equation}
\nabla _{X}\xi =-\varphi AX,  \label{eq-del-xi}
\end{equation}
\end{prop}

\noindent {\bf Proof.} Using (\ref{eq-JX}), (\ref{eq-JN=xi}), (\ref{eq-GF})
and (\ref{eq-WF}) in $(\tilde{\nabla}_{X}J)Y=0$, we get
\[
0=\left( \nabla _{X}\varphi \right) Y-\eta \left( Y\right) AX-h(X,Y)\xi
+\left( \left( \nabla _{X}\eta \right) Y\right) N+h(X,\varphi Y)N
\]%
Equating tangential and normal parts we get (\ref{eq-del-phi}) and (\ref%
{eq-del-eta}), respectively. Eq.~(\ref{eq-del-eta}) implies (\ref{eq-del-xi}%
). $\blacksquare $ \medskip

Now we obtain the following theorem of characterization for $\left(
\varepsilon \right) $-para Sasakian hypersurfaces.

\begin{th}
Let $M$ be an orientable hypersurface of an indefinite locally Riemannian
product manifold. Then $M$ is an $\left( \varepsilon \right) $-para Sasakian
manifold if and only if the shape operator is given by
\begin{equation}
A=-\varepsilon I+\varepsilon \eta \otimes \xi .  \label{eq-A-value}
\end{equation}
\end{th}

\noindent {\bf Proof.} Let $M$ be an $\left( \varepsilon \right) $-para
Sasakian manifold. By using (\ref{eq-s-pcm-def}) and (\ref{eq-del-xi}) we
get
\begin{equation}
AX=-\varepsilon X+\varepsilon \eta \left( X\right) \xi +\eta \left(
AX\right) \xi  \label{eq-AX}
\end{equation}%
In particular, we have $A\xi =\eta \left( A\xi \right) \xi $. Thus, we have
\begin{equation}
\eta \left( AX\right) =\varepsilon g\left( \xi ,AX\right) =\varepsilon
g\left( A\xi ,X\right) =\varepsilon g\left( \eta \left( A\xi \right) \xi
,X\right) =\eta \left( A\xi \right) \eta \left( X\right) .
\label{eq-eta(AX)-1}
\end{equation}%
Using this in (\ref{eq-AX}) we get
\begin{equation}
A=-\varepsilon I+\left( \varepsilon +\eta \left( A\xi \right) \right) \eta
\otimes \xi .  \label{eq-A-1}
\end{equation}%
Now, we use (\ref{eq-A-1}) in (\ref{eq-del-phi}) to find
\begin{equation}
\left( \nabla _{X}\varphi \right) Y=-\varepsilon \eta (Y)X+2\varepsilon \eta
\left( X\right) \eta (Y)\xi +2\eta \left( A\xi \right) \eta \left( X\right)
\eta (Y)\xi -g\left( X,Y\right) \xi .  \label{eq-del-phi-1}
\end{equation}%
From (\ref{eq-del-phi-1}) and (\ref{eq-eps-PS-def-1}) we get $\eta \left(
A\xi \right) =0$, which when used in (\ref{eq-A-1}) yields (\ref{eq-A-value}%
).

Conversely, using (\ref{eq-A-value}) in (\ref{eq-del-phi}) we see that $M$
is $\left( \varepsilon \right) $-para Sasakian manifold. $\blacksquare $
\medskip

Now, assume that the indefinite almost product Riemannian manifold $\tilde{M}
$\ is of almost constant curvature \cite{Yano-65-book} so that its curvature
tensor $\tilde{R}$ is given by
\begin{eqnarray}
\tilde{R}(X,Y,Z,W) &=&k\left\{ \tilde{g}\left( Y,Z\right) \tilde{g}\left(
X,W\right) -\tilde{g}\left( X,Z\right) \tilde{g}\left( Y,W\right) \right.
\nonumber \\
&&\qquad +\tilde{g}(JY,Z)\tilde{g}(JX,W)-\tilde{g}(JX,Z)\tilde{g}(JY,W)\}
\label{eq-ACCurv}
\end{eqnarray}%
for all vector fields $X,Y,Z,W$ on $\tilde{M}$. If $M$ is an $\left(
\varepsilon \right) $-para Sasakian hypersurface, then in view of (\ref%
{eq-A-value}) and (\ref{eq-ACCurv}) the Gauss equation becomes
\begin{eqnarray}
R(X,Y,Z,W) &=&\left( k-1\right) \left\{ g\left( Y,Z\right) g\left(
X,W\right) -g\left( X,Z\right) g\left( Y,W\right) \right\}   \nonumber \\
&&+k\left\{ g(\varphi Y,Z)g(\varphi X,W)-g(\varphi X,Z)g(\varphi
Y,W)\right\}   \nonumber \\
&&-\varepsilon g\left( Y,Z\right) \eta \left( X\right) \eta (W)-\varepsilon
g\left( X,W\right) \eta \left( Y\right) \eta (Z)  \nonumber \\
&&+\varepsilon g\left( X,Z\right) \eta \left( Y\right) \eta (W)+\varepsilon
g\left( Y,W\right) \eta \left( X\right) \eta (Z).  \label{eq-Gauss-1}
\end{eqnarray}%
After calculating $R\left( X,Y\right) \xi $ from (\ref{eq-Gauss-1}) and
comparing the resulting expression with \cite{TKYK-09-IJMMS}%
\begin{equation}
R\left( X,Y\right) \xi =\eta \left( X\right) Y-\eta \left( Y\right) X,
\label{eq-eps-PS-R(X,Y)xi}
\end{equation}%
we find that $k=2-\varepsilon $. With this value of $k$, from (\ref%
{eq-Gauss-1}), we obtain
\[
S=\left( \left( 2-\varepsilon \right) \left( n-2\right) -n\right) g+\left(
2-\varepsilon \right) \,{\rm trace}(\varphi )\,\Phi +\varepsilon
\,(4-\varepsilon -n)\eta \otimes \eta .
\]

Thus we have proved the following:

\begin{th}
An $(\varepsilon )$-para Sasakian hypersurface of an indefinite locally
Riemannian product manifold of almost constant curvature $\left(
2-\varepsilon \right) $ is Einstein like.
\end{th}

\begin{rem-new}
A hypersurface is called a {\em quasi-umbilical hypersurface} \cite{Chen-73}
if
\[
h\left( X,Y\right) =\alpha \;g\left( X,Y\right) +\beta \;u\left( X\right)
u\left( Y\right) ,
\]%
where $\alpha $ and $\beta $ are some smooth functions and $u$ is a $1$%
-form. From (\ref{eq-A-value}) we see that the $\left( \varepsilon \right) $%
-para Sasakian hypersurface is quasi-umbilical.
\end{rem-new}

\noindent Sadik Kele\c{s}

\noindent Department of Mathematics, Faculty of Arts and Sciences, \.{I}n%
\"{o}n\"{u} University

\noindent 44280 Malatya, Turkey

\noindent Email: skeles@inonu.edu.tr

\smallskip

\noindent Erol Kili\c{c}

\noindent Department of Mathematics, Faculty of Arts and Sciences, \.{I}n%
\"{o}n\"{u} University

\noindent 44280 Malatya, Turkey

\noindent Email: ekilic@inonu.edu.tr

\smallskip

\noindent Mukut Mani Tripathi

\noindent Department of Mathematics, Banaras Hindu University

\noindent Varanasi 221005, India.

\noindent Email: mmtripathi66@yahoo.com

\smallskip

\noindent Selcen Y\"{u}ksel Perkta\c{s}

\noindent Department of Mathematics, Faculty of Arts and Sciences, \.{I}n%
\"{o}n\"{u} University

\noindent 44280 Malatya, Turkey

\noindent Email: selcenyuksel@inonu.edu.tr

\end{document}